%hss.tex: Integrals Involving Rudin-Shapiro Polynomials and Sketch of a Proof of Saffari's Conjecture
%by Shalosh B. Ekhad and Doron Zeilberger
%Plain TeX
%begin macros

\baselineskip=14pt
\parskip=10pt

\font\eightrm=cmr8 

\magnification=\magstephalf

\def\1{{\overline{1}}}
\def\2{{\overline{2}}}
\parindent=0pt
\overfullrule=0in

\def\frac#1#2{{#1 \over #2}}
%\headline={\rm  \ifodd\pageno  \RightHead  \else  \LeftHead  \fi}
%\def\RightHead{\centerline{
%Title
%}}
%\def\LeftHead{ \centerline{Doron Zeilberger}}
%end macros

\bf
\centerline
{
Integrals Involving Rudin-Shapiro Polynomials and Sketch of a Proof of Saffari's Conjecture
}
\rm
\bigskip
\centerline
{\it By Shalosh B. EKHAD and Doron ZEILBERGER}
\bigskip
\qquad 

{\it \qquad \qquad
Dedicated to Krishnaswami ``Krishna'' Alladi, the tireless apostle of Srinivasa Ramanujan, yet a great mathematician on his own right.}

{\bf Preface: Krishna Alladi}

One of the greatest {\it disciples} of  Srinivasa Ramanujan, who did so much to make him a household name
in the mathematical community, and far beyond, is Krishnaswami ``Krishna'' Alladi.
Among many other things, he founded and is still editor-in-chief, of the  very successful {\it Ramanujan Journal}
(very ably managed by managing editor Frank Garvan), and initiated the SASTRA Ramanujan prize given to
promising young mathematicians.

But Krishna is not {\it just}  a mathematical leader, he is also a great number-theorist with very broad interests,
including analytic number theory and, inspired by Ramanujan,  $q$-series and partitions. That's why it is
not surprising that the conference to celebrate his 60th birthday, that took place last March, attracted attendees and
speakers with very diverse interests, and enabled the participants to learn new things
far afield from their own narrow specialty. That's how we found out, and got hooked on, {\it Rudin-Shapiro polynomials}.

{\bf Hugh Montgomery's Erd\H{o}s's colloquium}

One of the highlights of the conference was a fascinating talk by the eminent Michigan number theorist
(and Krishna's former postdoc mentor) Hugh Montgomery, who talked about {\it Littlewood polynomials} of
interest {\bf both} in pure number theory and, surprisingly, in {\it signal processing}. These are polynomials
whose coefficients are in $\{-1,1\}$. Among these stand out the famous {\it Rudin-Shapiro} polynomials,
introduced ([S1][S2]) by Harold ``Silent'' Shapiro\footnote{$^{1}$}
{\eightrm Harold S. Shapiro (S. originally stood for Seymour) was one of a brilliant cohort of students at City College,
in the late 1940s, that included Leon Ehrenpreis, Donald Newman, Israel Aumann, and another Harold Shapiro,
Harold N. Shapiro (N. originally stood for Nathaniel). But their friends, in order to distinguish between the two
Harold Shapiros, called them ``Silent'' and ``Noisy'' respectively. It is ironic that Harold Silent Shapiro's son
is the eminent, {\bf but very loud}, MIT cosmologist, Max Tegmark.
} and rediscovered by Walter Rudin ([R]).

{\bf The Rudin-Shapiro polynomials}

The Rudin-Shapiro polynomials, $P_k(z)$,  are best defined by the functional recurrence (see [Wi])
$$
P_k(z)=P_{k-1}(z^2) + zP_{k-1} (-z^2) \quad,
\eqno(DefiningRecurrence)
$$
with the initial condition $P_0(z)=1$.

As Hugh Montgomery described so well in his talk, these have amazing properties.
Both number-theorists and signal-processors are very interested in the so-called 
{\it sequence of (even) moments}, whose definition
usually involves the integral sign, but is better phrased entirely in terms of high-school algebra as follows.
$$
M_n(k):=CT [ P_k(z)^n P_k(z^{-1})^n ] \quad,
$$
where $CT$ denotes the ``constant term functional'', that for any Laurent polynomial $f(z)$ of $z$, extracts the
coefficient of $z^0$. For example $CT [4/z^2 +11/z+ 101 +5z +11z^{15}]=101$.

Can we find closed-form expressions for $M_n(k)$, in $k$, for any given, specific, positive integer $n$? Failing this, can we find
explicit expressions for the generating functions
$$
R_n(t) := \sum_{k=0}^{\infty} M_n(k) t^k \quad ?
$$

The sequence $M_1(k)$ has a very nice closed-form, $M_1(k)=2^k$. This is not very hard, even for humans.
Indeed, using Eq. $(DefiningRecurrence)$, we get
$$
P_k(z) P_k( z^{-1})= \left ( \, P_{k-1}(z^2) + zP_{k-1} (-z^2) \, \right ) \cdot  \left ( \, P_{k-1}(z^{-2}) + z^{-1} P_{k-1} (-z^{-2} ) \,
\right )
$$
$$
=P_{k-1}(z^2)P_{k-1}(z^{-2}) \, + \, P_{k-1}(-z^2)P_{k-1}(-z^{-2}) \, + \, 
\{ zP_{k-1}(-z^2)P_{k-1}(z^{-2}) + z^{-1} P_{k-1}(z^2)P_{k-1}(-z^{-2}) \} \quad .
$$
The quantity in the braces only has {\bf odd} powers, so its constant term vanishes. Hence
$$
M_1(k)=CT \, [ \, P_k(z) P_k(z^{-1}) \, ] \, = \, CT \, [ \, P_{k-1}(z^2)P_{k-1}(z^{-2}) \, ] 
\, + \, CT\, [\, P_{k-1}(-z^2)P_{k-1}(-z^{-2})\, ] \quad.
$$
Replacing $z^2$ by $z$ in the first term on the right, and $-z^2$ by $z$ in the second term, does not
change the constant term, hence,  we have the {\bf linear recurrence equation with constant coefficients}
$$
M_1(k)=2M_1(k-1) \quad,
$$
with the obvious initial condition $M_1(0)=1$, that implies the explicit expression $M_1(k)=2^k$.
Equivalently, the generating function $R_1(t)$ is given by
$$
R_1(t)=\frac{1}{1-2t} \quad .
$$

Let's move on to find an explicit formula for $M_2(k)$ and/or $R_2(t)$.
That was already done by smart human John Littlewood ([L]) but let's do it again.

Once again, let's use the defining recurrence for the Rudin-Shapiro polynomials, but let's abbreviate
$$
a(k)(z)=P_{k}(z) \quad , \quad   b(k)(z)=P_{k}(-z) \quad , \quad  A(k)(z)=P_{k}(z^{-1}) \quad , \quad  B(k)(z)=P_{k}(-z^{-1}) \quad .
$$

We have
$$
P_k(z)^2 P_k( z^{-1})^2= \left ( \, P_{k-1}(z^2) + zP_{k-1} (-z^2) \, \right )^2 
\cdot  
\left ( \, P_{k-1}(z^{-2}) + z^{-1} P_{k-1} (-z^{-2}) \, \right )^2 \quad .
$$

Expanding, discarding odd terms, replacing $z^2$ by $z$,
and using trivial symmetries due to the fact that the functional  $CT$ is preserved
under the dihedral group $\{ z \rightarrow z,  z \rightarrow -z,  z \rightarrow z^{-1},  z \rightarrow -z^{-1}\}$, we get that
$$
CT \, [ \, a(k)^2 A(k)^2]=2 CT[ a(k-1)^2 A(k-1)^2 \, ] \,  -2 CT \, [ \,  z a(k-1)^2 B(k-1)^2] + 4 CT [ a(k-1)A(k-1)b(k-1)B(k-1) \, ] \quad .
$$

The first term is an old friend, our quantity of interest with $k$ replaced by $k-1$, but the other two are newcomers.
So we do the same treatment to them. They in turn, may (and often do) introduce new quantities, but if all goes well,
there would only be finitely many sequences, and we would get a {\bf finite} system of first-order linear recurrences.
This indeed happens, and one gets, for the generating functions of the encountered sequences, 
a system of six equations with six unknowns,
and in particular, we get (in a split second, of course, we let Maple do it) that our desired object,
the generating function of the squence $CT[a(k)^2A(k)^2]$, alias, $R_2(t)$,  is given by:
$$
R_2(t)={\frac {4\,t \, + \, 1}{ \left( 1\,+ \, 2\,t \right)  \left( 1 \, - \, 4\,t \right) }} \quad
=
\frac{4}{3}\, \frac{1}{1-4t}- \frac{1}{3} \, \frac{1}{1 \, + \, 2\,t} \quad .
$$
By extracting the coefficient of $t^k$,  we even get a nice explicit expression for $M_2(k)$,
already known to Littlewood
$$
M_2(k) = \frac{4}{3} 4^k  -\frac{1}{3} (-2)^k \quad .
$$

This can be done for {\it any} monomial 
$$
z^{\alpha_0} a(k)^{\alpha_1} A(k)^{\alpha_2} b(k)^{\alpha_3} B(k)^{\alpha_4} \quad.
$$
Define the sequence
$$
E[\alpha_0, \alpha_1, \alpha_2, \alpha_3, \alpha_4](k):=
CT \, [ \, z^{\alpha_0} a(k)^{\alpha_1} A(k)^{\alpha_2} b(k)^{\alpha_3} B(k)^{\alpha_4} \, ] \quad .
$$
Replacing $a(k), A(k), b(k), B(k)$ by their expressions in terms of $z,a(k-1), A(k-1), b(k-1), B(k-1)$, expanding,
discarding odd terms, replacing $z^2$ by $z$, and replacing each monomial by its {\it canonical form},
implied by the above-mentioned action of the dihedral group that preserves CT, we can express, each such $E[.]$, 
in terms of other 
$E[.]$'s evaluated at $k-1$. It is possible to show 
(and it has been done by Doche and Habsieger [DH], using a different approach) that this process terminates,
and eventually we will not get any new sequences, leaving us with a finite system of linear equations
for the corresponding generating functions, that can be automatically solved, and lead to an expression
in terms of a rational function, since we get a first-order system
$$
{\bf F}(t)= {\bf v} + t {\bf A} {\bf F}(t) \quad ,
$$
(where ${\bf F}(t)$ is the vector of generating functions whose first component is our desired one),
for some matrix ${\bf A}$, of integers that the computer finds automatically, and our object of desire is
the first component of ${\bf F}(t)=(1-t{\bf A})^{-1} {\bf v}$.

While it is painful for a human to do this, a computer does not mind, and the Maple package

{\tt HaroldSilentShapiro.txt}

accompanying this article does it for any desired monomial in $z,P_k(z), P_k(-z),  P_k(z^{-1}),  P_k(-z^{-1})$.
See the output files accompanying this article, that may be viewed from the front of this article

{ \tt http://www.math.rutgers.edu/\~{}zeilberg/mamarim/mamarimhtml/hss.html} \quad .

Unlike the beautiful approach of Doche and Habsieger, that uses clever human pre-processing to establish
an {\it algorithm}, that was then hard-programmed by hand, our approach is naive ``dynamical programming'',
where we don't make any {\it a priori} human analysis, and let the computer introduce new quantities as needed.
To guarantee that it {\it halts}, we input a parameter, that we call {\tt K},  and if the size of the system
exceeds $K$ it returns {\tt FAIL}, leaving us the option to forget about it, or try again with a larger {\tt K}.

{\bf Higher moments and Saffari's Conjecture}

Now that we have  reduced, for any specific positive integer $n$, the computation of the generating function of the sequence of  moments $M_n(k)$,
that we call $R_n(t)$, to a routine calculation, we can ask our beloved computer to crank-out as many
of them as it can output in a reasonable amount of time. According to the output file

{\tt http://www.math.rutgers.edu/\~{}zeilberg/tokhniot/oHaroldSilentShapiro1.txt} \quad ,

we get 

$$
R_1(t)= \frac{1}{1-2t} \quad ,
$$
$$
R_2(t) ={\frac {1+4\,t}{ \left( 1 \, + \, 2\,t \right)  \left( 1 \, - \, 4 \,t \right) }} \quad ,
$$
[both of which were already given above],
$$
R_3(t)= {\frac {1+16\,t}{ \left( 1+4\,t \right)  \left( 1- 8\,t \right) }} \quad ,
$$
$$
R_4(t)= 
$$
-(90194313216*t**11-15300820992*t**10-1979711488*t**9-292552704*t**8-
22216704*t**7+10649600*t**6-1024*t**5-144384*t**4+7008*t**3+664*t**2-54*t-1)/
((8*t+1)*
(16*t-1)*(1409286144*t**10-264241152*t**9-25690112*t**8-4128768*t**7-311296*t**6+
170496*t**5-2624*t**4-2208*t**3+148*t**2+8*t-1)),

\vfill\eject

$$
R_5(t)=
$$
-(369435906932736*t**11-32160715112448*t**10-2001454759936*t**9-  \hfill\break
145223581696*t**8-4454350848*t**7+1392508928*t**6-5865472*t**5-4599808*t**4  \hfill\break
+123648*t**3+4768*t**2-220*t-1)/\hfill\break
((1+16*t)*(32*t-1)*(1443109011456*t**10-135291469824*t**9 
-6576668672*t**8-528482304*t**7\hfill\break
-19922944*t**6+5455872*t**5-41984*t**4-17664*t**3+592*t**2+16*t-1))  \quad .

To see $R_k(t)$ for $6 \leq k \leq 10$, look at the above-mentioned output file. Of course, one can easily go further.
Note that these have already been computed in [DH] (but their output  is not easily accessible to the casual reader).

By looking at the smallest root of the  denominator of $R_k(t)$ and computing the residue, one 
confirms for small (and not so small!) values of $k$ (and one can easily go much further), the following conjecture of
Bahman Saffari, as already done in [DH] (for small $k$).

{\bf Saffari's Conjecture} 
$$
M_n(k) \, \sim \, \frac{2^n}{n+1} \cdot (2^n)^k \quad .
$$

Saffari never published his conjecture, and it is mentioned as ``private communication'' in [DH].

{\bf Sketch of a proof of Saffari's Conjecture}

While for each {\it numeric} $n$, one can get  an explicit expression, in {\it symbolic} $t$, for $R_n(t)$, 
these get more and more complicated as $n$ gets larger, and there is (probably) no hope to get
an explicit expression, in {\bf symbolic} $n$, for $R_n(t)$, from which one can deduce that the smallest root
(in absolute value) of the denominator is $2^{-n}$ and the residue is $\frac{2^{n}}{n+1}$.

But one can prove {\it rigorously} Saffari's conjectured asymptotic formula as follows.

Let $n$ be a general (symbolic) positive integer. Recall that we are interested in the sequence 
$$
M_n(k):=CT [P_k(z)^n P_k(z^{-1})^n] \quad ,
$$
that we abbreviate
$$
a^n A^n \quad,
$$
under the convention 
$$
a=P_k(z) \quad , \quad  b=P_k(-z) \quad , \quad A=P_k(z^{-1}) \quad , \quad  B=P_k(-z^{-1}) \quad.
$$
To get a scheme we use the {\it rewriting rules}, implied by the defining recurrence
$$
a \rightarrow a+zb \quad , \quad
b \rightarrow a-zb \quad , \quad
A \rightarrow A+z^{-1}B \quad , \quad
B \rightarrow A -z^{-1}B \quad , 
$$
where the discrete argument on the left is $k$ and on the right  $k-1$, and the continuous argument on the left is $z$
and on the right is $z^2$.

Using the binomial theorem, we have
$$
a^n A^n \rightarrow (a+zb)^n (A+ z^{-1} B)^n =
\left ( \sum_{i=0}^{n} {{n} \choose {i}} a^i (zb)^{n-i} \right )
\left ( \sum_{j=0}^{n} {{n} \choose {j}} A^j (z^{-1}B)^{n-j} \right )
$$
$$
=\sum_{i=0}^{n} \sum_{j=0}^{n} {{n} \choose {i}} {{n} \choose {j}} a^i b^{n-i} A^j B^{n-j} z^{j-i} 
$$
$$
=\sum_{i=0}^{n} {{n} \choose {i}}^2 (aA)^i (bB)^{n-i} \quad + SmallChange \quad ,
$$
where $SmallChange$ is a linear combination of unimportant monomials and 
we {\bf define} an {\bf important monomial} (in $a,A,b,B,z$) to be any member of the set of monomials
$$
\{ (aA)^m (bB)^{n-m} \, \vert \, 0 \leq m \leq n \} \quad .
$$

Let's try to find the ``going down'' evolution-step for the other important monomials.

We have
$$
(aA)^m (bB)^{n-m} \rightarrow (a+zb)^m  (A+ z^{-1}B)^m  (a-zb)^{n-m}  (A- z^{-1}B)^{n-m}
$$
$$ 
= \left ( \sum_{i_1=0}^{m} {{m} \choose {i_1}} a^{i_1} (zb)^{m-i_1} \right )
\left ( \sum_{i_2=0}^{m} {{m} \choose {i_2}} A^{i_2} (z^{-1}B)^{m-i_2} \right ) \cdot
$$
$$
\left ( \sum_{i_3=0}^{n-m} {{n-m} \choose {i_3}} a^{i_3} (-zb)^{n-m-i_3} \right )
\left ( \sum_{i_4=0}^{n-m} {{n-m} \choose {i_4}} A^{i_4} (-z^{-1}B)^{n-m-i_4} \right )
$$
$$
= \, \sum_{i_1=0}^{m} \sum_{i_2=0}^{m} \sum_{i_3=0}^{n-m} \sum_{i_4=0}^{n-m}
 {{m} \choose {i_1}} {{m} \choose {i_2}} {{n-m} \choose {i_3}} {{n-m} \choose {i_4}}
(-1)^{i_3+i_4} a^{i_1+i_3} A^{i_2+i_4}b^{n-i_1-i_3}B^{n-i_2-i_4} z^{i_2 - i_1 +i_4 - i_3} \quad .
$$

The coefficient of a typical important monomial, $(aA)^r (bB)^{n-r}$ ($0 \leq r \leq n$) in the above quadruple sum is
$$
\sum_{i_1+i_3=r \, , \, i_2+i_4=r}
(-1)^{i_3+i_4}
 {{m} \choose {i_1}} {{m} \choose {i_2}} {{n-m} \choose {i_3}} {{n-m} \choose {i_4}}
$$
$$
= \, \sum_{i_1=0}^{r}\sum_{i_2=0}^{r}
(-1)^{i_1+i_2}   {{m} \choose {i_1}} {{m} \choose {i_2}} {{n-m} \choose {r-i_1}} {{n-m} \choose {r-i_2}}
$$
$$
= \, \left ( \sum_{i_1=0}^{r} (-1)^{i_1}  {{m} \choose {i_1}}  {{n-m} \choose {r-i_1}} \right )
\left ( \sum_{i_2=0}^{r} (-1)^{i_2}  {{m} \choose {i_2}}  {{n-m} \choose {r-i_2}} \right )
$$
$$
=\left ( \sum_{i=0}^{r} (-1)^{i}  {{m} \choose {i}}  {{n-m} \choose {r-i}} \right )^2 \quad .
$$

This is an important quantity, so let's give it a name
$$
K_n(m,r):=\left ( \sum_{i=0}^{r} (-1)^{i}  {{m} \choose {i}}  {{n-m} \choose {r-i}} \right )^2 \quad .
$$

All the remaining monomials belong to $SmallChange$, and we have the  general `evolution equation'
$$
(aA)^m (bB)^{n-m} \rightarrow \,
\sum_{r=0}^{n} K_n(m,r) (aA)^r (bB)^{n-r} + \quad SmallChange \quad.
$$

Assuming for now, that  $SmallChange$ is, asymptotically less than the ``important monomials''
(i.e. the rate of growth of a small-change sequence divided by an ``important monomial'' sequence is $o(1)$), let
$\alpha_n$ be the largest eigenvalue of the $n+1$ by $n+1$ matrix $K_n$ (whose $(m,r)$ entry is $K_n(m,r)$), then
for $0 \leq m \leq n$
$$
CT[ (P_k(z)P_k(z^{-1}))^m (P_k(-z)P_k(-z^{-1}))^{n-m} ] \, \sim  \, c_m (\alpha_n)^k \quad,
$$
where $(c_0, \dots ,c_n)$ is an eigenvector corresponding to the largest eigenvalue, $\alpha_n$.

We now need two elementary propositions that should be provable  using
the {\bf Wilf-Zeilberger algorithmic proof theory}, ([Z][KKZ] for the first, [WZ][AZ] for the second).
They may be even provable by purely human means, but since we know {\it for sure} that they are
both true, we do not bother.

{\bf Proposition 1.} The characteristic polynomial, $\det(z {\bf I}- K_n)$, of the $(n+1) \times (n+1)$ matrix $K_n$ whose $(m,r)$ entry is
$$
K_n(m,r):=\left ( \sum_{i=0}^{r} (-1)^{i}  {{m} \choose {i}}  {{n-m} \choose {r-i}} \right )^2 \quad ,
$$
equals
$$
\det(z{\bf I} -K_n)=
z^{\lfloor (n+1)/2 \rfloor}
\prod_{j=0}^{\lfloor n/4 \rfloor} 
\left ( \, z- 2^{n-4j} {{4j} \choose {2j}} \, \right )
\prod_{j=0}^{\lfloor (n-2)/4 \rfloor} \left ( \,  z+ 2^{n-4j-2} {{4j+2} \choose {2j+1}} \, \right ) \quad .
$$

[To confirm this {\it shaloshable} determinant identity for $n \leq N$, type, in the Maple package \hfill\break
{\tt HaroldSilentShapiro.txt}, {\tt CheckCP(N); } .
For example, {\tt CheckCP(20);} returns {\tt true} in one second, and  {\tt CheckCP(40);} returns {\tt true} in $20$ seconds.]

So the non-zero eigenvalues of the matrix $K_n$ are
$$
\{\,\,  2^{n-4j} {{4j} \choose {2j}} \, \, ;  \, \, 0 \leq j \leq \lfloor n/4 \rfloor \,\, \} \, \bigcup \,
\{ \,\,  -2^{n-4j-2} {{4j+2} \choose {2j+1}} \, \, ;  \, \, 0 \leq j \leq \lfloor (n-2)/4 \rfloor \,\, \} \quad .
$$

In particular, the largest eigenvalue (in absolute value) is indeed $2^n$. We also need the following
{\it shaloshable} binomial coefficients identity.

{\bf Proposition 2.} The vector $(c_0, \dots, c_n)$ defined by $c_r={{n} \choose {r}}^{-1}$ ($0 \leq r \leq n$)
is an eigenvector of the matrix $K_n$ corresponding to its largest eigenvalue $2^n$. In other words,
for $0 \leq m \leq n$
$$
\sum_{r=0}^{n} K_n(m,r) c_r =2^n c_m \quad .
$$

[To confirm this {\it shaloshable} binomial coefficient identity for $n \leq N$, type, in the Maple package
{\tt HaroldSilentShapiro.txt}, {\tt CheckEV(N); }.
For example, {\tt CheckEV(50);} returns {\tt true} in two seconds, and  {\tt CheckCP(100);} returns {\tt true} in $30$ seconds.]

But an eigenvector is only determined up to a constant multiple. Let's find it (modulo the Small Change hypothesis).
We know that
$$
CT \, [ \, (P_k(z)P_k(z^{-1}))^m (P_k(-z)P_k(-z^{-1}))^{n-m} \,] \,  \sim \, \frac {C}{ { {n} \choose {m}}}  \cdot (2^n)^k \quad ,
$$
for {\it some} constant $C$.  To find it, we use the well-known, and easily proved identity (see [Wi])
$$
P_k(z)P_k(z^{-1}) + P_k(-z)P_k(-z^{-1}) =2^{k+1} \quad.
$$
Raising it to the $n$-th power, using the binomial theorem, and taking the constant term, we have
$$
\sum_{m=0}^{n} {{n} \choose {m}} CT[(P_k(z)P_k(z^{-1}))^m ( P_k(-z)P_k(-z^{-1}))^{n-m}] =2^{(k+1)n} \quad.
$$

Hence
$$
\sum_{m=0}^{n} {{n} \choose {m}} \frac {C}{ { {n} \choose {m}}} \cdot (2^n)^k \quad  =2^{(k+1)n} \quad,
$$
that implies that
$$
C=\frac{2^n}{n+1} \quad .
$$

We just established

{\bf Proposition 3}: Modulo the Small Change Hypothesis, for $0 \leq m \leq n < \infty$
$$
CT[(P_k(z)P_k(z^{-1}))^m (P_k(-z)P_k(-z^{-1}))^{n-m} \, ] \, \sim \, \frac {2^n} { (n+1){ {n} \choose {m}}} \cdot (2^n)^k \quad .
$$
In particular, taking $m=n$, we get Saffari's conjecture (for even moments)
$$
M_n(k)=CT[P_k(z)^n P_k(z^{-1})^n] \, \sim \, \frac {2^n}{n+1} \cdot (2^n)^k \quad .
$$

{\bf Towards a Proof of the Small Change Hypothesis}

It would have been great if the ``children'' of each unimportant monomial, in the evolution equation
described above (implemented in procedure {\tt GD} in our Maple package), would all be
unimportant. Then we could have easily proved, by induction that, asymptotically, they are insignificant
compared to the important monomials. It turns out that for {\it most} unimportant monomials, this is indeed
the case, but there are a few, that we call {\it false pretenders} that do have important children.

It should not be hard to fully characterize these. In fact it turns out (empirically, for now) that for $n$ even there
are $(n/2)^2-1$ of then, and for $n$ odd there are $(n^2-1)/4$. Then for those false pretenders one should
be able to describe all their important children, and then prove that the leading terms of their
contributions cancel out (using the inductive hypothesis, and Prop. 3).

This has been verified empirically up to $n \leq 16$. See procedures {\tt Medio} and {\tt MedioP} in the Maple
package  {\tt HaroldSilentShapiro.txt}  \quad .

{\bf Hugh Montgomery's Stronger Conjecture}

In [M], Hugh Montgomery considered the more general sequences
$$
M_{m,n}(k) := CT \, [ \,  P_k(z)^m P_k(z^{-1})^n \, ] \quad . 
$$
He conjectured that, for $m \neq n$ ,
$$
M_{m,n}(k) = o (\, 2^{(m+n)k/2} \, ) \quad .
$$
Once again, the generating function, for each specific $m$ and $n$, is always a rational function, and our
Maple package (procedure {\tt RS(m,n,t,K)}) computes them, and procedure {\tt MamarH(N,K,t)} 
prints out an article confirming Hugh Montgomery's conjecture, as well as giving the generating
functions for $1 \leq m < n \leq N$. ({\tt K} is a parameter that  should be made large enough, say 1000).

To see the output for  $1 \leq m < n \leq 7$, go to:\hfill\break
{\tt http://www.math.rutgers.edu/\~{}zeilberg/tokhniot/oHaroldSilentShapiro2.txt} \quad ,\hfill\break
that contains the explicit expressions for all these cases, and confirms Montgomery's conjecture with
a vengeance. Unlike the $m=n$ case, the smallest root (alias the reciprocal of the largest eigenvalue) is not `nice', and
there are usually several roots with smallest  absolute value, hence the sequences
often oscillate. Nevertheless, Montgomery's conjecture is true for all  $1 \leq m < n \leq 7$, 
and one could go much further.

{\bf Let's Generalize!}

The same approach works for {\it any} sequence of Laurent polynomials defined by a recurrence of the form
$$
P_k(z)=C_1(z) P_{k-1}(z^r) + C_2(z) P_{k-1} (-z^r)+C_3(z) P_{k-1}(z^{-r}) + C_4(z) P_{k-1} (-z^{-r}) \quad,
$$
with the initial condition $P_0(z)=1$, where $C_1(z), C_2(z), C_3(z), C_4(z)$ are Laurent polynomials of
degree less than $r$ and low-degree larger than $-r$, for {\it any} positive integer $r$ larger than $1$.

One always gets a finite scheme, and hence a rational generating function for the sequence
$$
S[\alpha_0, \alpha_1, \alpha_2, \alpha_3, \alpha_4](k):=
CT \left [ \, z^{\alpha_0} P_k(z)^{\alpha_1} P_k(z^{-1})^{\alpha_2} P_k(-z)^{\alpha_3} P_k(-z^{-1})^{\alpha_4} \right ] \quad ,
$$
for any non-negative $\alpha_0, \alpha_1, \alpha_2, \alpha_3, \alpha_4$. This is implemented in the
Maple package \hfill\break
{\tt ShapiroGeneral.txt} also available from the webpage  of this article, or directly from \hfill\break
{\tt http://www.math.rutgers.edu/\~{}zeilberg/tokhniot/ShapiroGeneral.txt } \quad .

{\bf Let's (not!) Generalize Even More!}

The set $\{1,-1\}$ is a multiplicative subgroup of the field of complex numbers.
For any finite multiplicative subgroup $G$ of the  field of complex numbers, and any positive integer $r$ larger than $1$,
the same approach should
be able to handle sequences of polynomials given by a recurrence
$$
P_k(z)= \sum_{g \in G} \alpha_g(z) P_{k-1}(g z^r) + \sum_{g \in G} \beta_g(z) P_{k-1}(gz^{-r})
 \quad , \quad P_0(z)=1 \quad ,
$$
where $\alpha_g(z),\beta_g(z)$ are $2|G|$ given Laurent polynomials in $z$ of degree $<r$ and low-degree $>-r$.

This includes the case treated in [D],  where $G$ is a cyclotomic group.

We could go even further, with {\it higher order} recurrences (as opposed to only first order),
several continuous variables (as opposed to only $z$), and, presumably, even several discrete variables
(as opposed to only $k$), but {\it enough is enough!}

{\bf Added May 27, 2016}: Brad Rodgers, independently, and simultaneously, found a (complete) proof of 
Saffari's conjecture, that he is writing up now and will soon post in the arxiv. 
Meanwhile, you can read his proof in a letter posted out in  \hfill\break
{\tt http://www.math.rutgers.edu/\~{}zeilberg/mamarim/mamarimhtml/BradleyRodgersLetter.pdf}.

{\bf Added June 7, 2016}: Brad Rodgers' beautiful paper, that also proves the more general Montgomery conjecture,
mentioned above, is now available here: {\tt http://arxiv.org/abs/1606.01637}.

{\bf References}

[AZ] Moa Apagodu and Doron Zeilberger,
{\it Multi-Variable Zeilberger and Almkvist-Zeilberger Algorithms and the Sharpening of Wilf-Zeilberger Theory},
Adv. Appl. Math. {\bf 37} (2006), 139-152. \hfill\break
{\tt http://www.math.rutgers.edu/\~{}zeilberg/mamarim/mamarimhtml/multiZ.html} \quad .

[D] Christophe Doche, {\it Even moments of generalized Rudin-Shapiro polynomials},
Math. Comp. {\bf 74}(2005), 1923-1935.  \hfill \break
{\tt http://www.ams.org/mcom/2005-74-252/S0025-5718-05-01736-9/S0025-5718-05-01736-9.pdf}  \quad [viewed May 6, 2016].

[DH] Christophe  Doche and Laurent Habsieger, {\it Moments of the Rudin-Shapiro polynomials},
Journal of Fourier Analysis and Applications {\bf 10} (2004), 497-505. \hfill \break
{\tt http://web.science.mq.edu.au/\~{}doche/049.pdf} \quad [viewed May 6, 2016].

[KKZ] Christoph Koutschan, Manuel Kauers, and Doron Zeilberger,
{\it A Proof Of George Andrews' and David Robbins' q-TSPP Conjecture},
Proceedings of the National Academy of Science, {\bf 108}, no. 6 (Feb. 8, 2011), 2196-2199. \hfill\break
{\tt http://www.math.rutgers.edu/\~{}zeilberg/mamarim/mamarimhtml/qtsppRig.html} \quad .

[L] John E. Littlewood, {\it ``Some Problems in Real and Complex Analysis''}, D.C. Heath and Co. Raytheon Education Co.,
Lexington, Mass., 1968.

[M] Hugh Montgomery, {\it  Problem}, in: {\it ``Problem Session'', 
p. 3029, Mathematisches Forschungsinstitut Oberwolfach,
Report No. 51/2013,  DOI: 10.4171/OWR/2013/51, Analytic Number Theory, 
20 October - 26 October 2013},  3029-3030. \hfill \break
{\tt https://www.mfo.de/document/1343/OWR\_2013\_51.pdf} \quad [viewed May 6, 2016].

[R] Walter Rudin, {\it Some theorems on Fourier coefficients}, Proc. Amer. Math. Soc. {\bf 10} (1959), 855-859.

[S1] Harold S. Shapiro, {\it ``Extremal problems for polynomial and power series''}, M.Sc. thesis, MIT (1951).  \hfill\break
{\tt http://dspace.mit.edu/bitstream/handle/1721.1/12198/30502786-MIT.pdf?sequence=2} \quad [viewed May 11, 2016].

[S2] Harold S. Shapiro, {\it ``Extremal problems for polynomial and power series''}, Ph.D. thesis, MIT (1952).  \hfill\break
{\tt http://dspace.mit.edu/bitstream/handle/1721.1/12247/30752144-MIT.pdf?sequence=2}
\quad [viewed May 11, 2016].

[Wi] The Wikipedia Foundation, {\it ``Shapiro polynomials''}.  \hfill\break
{\tt https://en.wikipedia.org/wiki/Shapiro\_polynomials} \quad [viewed May 13, 2016].

[WZ] Herbert S. Wilf  and Doron Zeilberger,
{\it  An algorithmic proof theory for hypergeometric (ordinary and "q") multisum/integral identities},
Invent. Math. {\bf 108}(1992), 575-633.
\hfill\break
{\tt http://www.math.rutgers.edu/\~{}zeilberg/mamarim/mamarimhtml/multiwz.html} \quad .

[Z] Doron Zeilberger,
{\it The Holonomic ansatz II: automatic discovery(!) and proof(!!) of Holonomic determinant evaluations},
Annals of Combinatorics {\bf 11}(2007), 241-247. \hfill\break
{\tt http://www.math.rutgers.edu/\~{}zeilberg/mamarim/mamarimhtml/ansatzII.html} \quad .

\vfill\eject

\bigskip
\bigskip
\hrule
\bigskip
Shalosh B. Ekhad, c/o D. Zeilberger, Department of Mathematics, Rutgers University (New Brunswick), Hill Center-Busch Campus, 110 Frelinghuysen
Rd., Piscataway, NJ 08854-8019, USA. 
\bigskip
\hrule
\bigskip
Doron Zeilberger, Department of Mathematics, Rutgers University (New Brunswick), Hill Center-Busch Campus, 110 Frelinghuysen
Rd., Piscataway, NJ 08854-8019, USA. \hfill \break
zeilberg at math dot rutgers dot edu \quad ;  \quad {\tt http://www.math.rutgers.edu/\~{}zeilberg/} \quad .
\bigskip
\hrule
\bigskip
{\bf  First Written:  May 20, 2016} ; 
{\bf  This version: June 7, 2016} .

\end